\documentclass{article}

\newtheorem{theorem}{Theorem}
\newtheorem{lemma}[theorem]{Lemma}

\newtheorem{definition}[theorem]{Definition}
\newtheorem{example}[theorem]{Example}
\newtheorem{remark}[theorem]{Remark}

\newcommand{\examp}[1]
  {\begin{example} {\rm #1} \end{example}}

\def\QED{\quad\blackslug\lower 8.5pt\null}

\begin{document}

\begin{center} \LARGE{\textbf{On the Structure of Submanifolds}}

\vspace*{3mm}

\LARGE{\textbf{with Degenerate Gauss Maps}}

\vspace*{3mm}

{\large  Maks A. Akivis and Vladislav V. Goldberg}
\end{center}

\vspace*{5mm}

{\footnotesize{\it Abstract}.
An $n$-dimensional submanifold $X$ of a projective space
$P^N (\mathbf{C})$ is called  tangentially degenerate
if the rank of its Gauss mapping $\gamma: X \rightarrow G (n, N)$
satisfies $0 < \mbox{{\rm rank}} \; \gamma < n$.

The authors systematically study the geometry of tangentially
degenerate submanifolds of a projective space $P^N (\mathbf{C})$.
By means of the focal images, three  basic types of submanifolds
are discovered: cones, tangentially degenerate  hypersurfaces,
and torsal submanifolds. Moreover, for  tangentially degenerate
submanifolds, a structural theorem is proven. By this theorem,
tangentially degenerate  submanifolds that do not belong to one
of the basic types are foliated into submanifolds of basic types.
In the proof the authors introduce irreducible, reducible, and
completely reducible  tangentially degenerate  submanifolds. It
is found  that cones and tangentially degenerate  hypersurfaces
are irreducible, and torsal  submanifolds  are completely
reducible while all other tangentially degenerate  submanifolds
not belonging to basic types are reducible.

\vspace*{2mm}

{\it Keywords}: tangentially degenerate submanifold,
submanifold with degenerate Gauss mapping,
structure theorem.

{\it $2000$ Subject Classification}: 53A20
}

\setcounter{equation}{0}

\vspace*{5mm}

\setcounter{section}{-1}

\section{Introduction}
An $n$-dimensional submanifold $X$ of a projective space $P^N
(\mathbf{C})$ is called  {\em tangentially degenerate} if the rank of its
Gauss mapping $\gamma: X \rightarrow G (n, N)$ is less than $n, \;
r = \mbox{{\rm rank}} \; \gamma < n$.
Here $x \in X, \; \gamma (x) = T_x (X)$, and $T_x (X)$ is the
tangent subspace to $X$ at $x$ considered as an $n$-dimensional
projective space $P^n$.  The number  $r$
is also called the {\em rank} of $X, \; r  =  \mbox{{\rm rank}} \; X$.
The case $r = 0 $ is trivial one: it gives just an $n$-plane.
A submanifold  $X$ is called tangentially degenerate if $0 \leq r <
n$, and it is denoted by $V^n_r, \; X = V^n_r$. The submanifolds of
rank $r < n$ have been the object of numerous investigations because of their
analogy to developable surfaces in a three-dimensional space and
because of their significance in the theory of the curvature of submanifolds.

The tangentially degenerate submanifolds $X$ of rank $r < n$ were
first considered by \'{E}. Cartan [C 16] in connection with his
study of metric deformations of hypersurfaces, and in [C 19] in
connection with his study of manifolds of constant curvature.
In particular, Cartan proved that if $V^n_r$ is a tangentially
degenerate submanifold  of dimension $n$ and rank $r$ and  the
dimension $\rho$ of the  osculating subspace is
$\rho = n + \frac{1}{2} r (r + 1)$, then $V^n_r$ is a cone with
an $(n - r - 1)$-dimensional vertex.

It appeared that  tangentially degenerate submanifolds
are less rigid under the affine and projective deformations.
 Yanenko  [Ya 53]  studied them in connection with his study
 of metric deformations of submanifolds.
  Akivis [A 57, 62] studied them in multidimensional projective
space, considered their focal images (the locus of singular points
and the locus of singular hyperplanes) and applied the latter to
clarify the structure of  the tangentially degenerate
submanifolds. Savelyev  [Sa 57, 60] found a classification of
tangentially degenerate  submanifolds and described in detail
the tangentially degenerate submanifolds of rank 2.  Ryzhkov
[R 58] showed that a tangentially degenerate submanifolds $X$
of rank $r$ can be constructed by using the Peterson
transformation of  $r$-dimensional submanifolds, and in [R 60]
he proved  that such a construction is quite general. That is,
 by means of it, an arbitrary   tangentially degenerate
 submanifold $X$ of rank $r$ can be obtained (see also the
 survey paper [AR 64]).  In particular, Ryzhkov [R 60]
 generalized the above mentioned Cartan's result to the case
 $ n + 1 + \frac{1}{2}  r (r - 1) <\rho < n + \frac{1}{2} r (r + 1)$.
 Brauner [Br 38], Wu [Wu 95], and  Fischer and Wu [FW 95]
   studied such submanifolds in an  Euclidean $N$-space.

For a submanifold $V^n$ of a Riemannian space $V^N, \; n < N$,
 Chern and Kuiper [CK 52] introduced the notion of the index
 of relative nullity $\mu (x)$, where $x \in V^n$
 (see also [KN 69], p. 348). The  submanifolds $V^n$, for which
 $\mu(x)$ is constant and greater than 0 for all points $x \in V^n$,
 are called {\em strongly parabolic}. Akivis [A 87] proved
 that if a space $V^N$ admits a projective realization
 (this is always the case for the simply connected Riemannian spaces
  of constant  curvature, see [W 72], Ch. 2) and if the index
  $\mu (x)$ is constant on $V^n \subset
 V^N$, then the index $\mu (x)$ is connected with the rank
 of a  submanifold $V^n$ by the relation $\mu (x) = n - r$.
 This implies that the results of the papers [A 57, 62]
 as well as the results of the current paper can be applied to
 the study of strongly parabolic submanifolds of the Euclidean and
 non-Euclidean spaces. In particular, in [A 87] Akivis
 proved the existence of tangentially degenerate submanifolds
 in such spaces without singularities
 and constructed examples of such submanifolds.
 The main  results of papers indicated above can be found
 in Chapter 4 of the book [AG 93].
In the same paper [A 87], Akivis also proved that
a hypersurface in a four-dimensional Euclidean space
$\mathbf{R}^4$ considered by Sacksteder [S 60] is of rank 2
and without singularities, and that this hypersurface
is a particular case of a series of examples presented
in [A 87]. Later Akivis and Goldberg [AG 00] proved
that a similar example constructed (but not
published) by Bourgain and published in [Wu 95],
[I 98, 99a, 99b], and [WZ 99] coincides with
Sacksteder's example up to a coordinate transformation.

Griffiths and Harris [GH 79] (Section  2, pp. 383--393)
considered  tangentially degenerate submanifolds from the
point of view algebraic geometry. They used  the term
``submanifolds with degenerate Gauss mappings'' instead
 of the term ``tangentially degenerate  varieties''.
 They used this term to avoid a confusion with submanifolds with
 degenerate tangential varieties considered in Section
  5 of the same paper.

Following [GH 79], Landsberg [L 96] considered tangentially
degenerate  submanifolds. His recently published  notes
[L 99] are in some sense an  update to the   paper [GH 79].
Section  5 (pp. 47--50) of these notes is devoted
to tangentially degenerate  submanifolds.

 Griffiths and Harris [GH 79] asserted a structure
theorem for submanifolds with degenerate Gauss mappings,
that is, for the varieties $X = V^n_r$ such that
$\dim \gamma (X)<\dim X$. They asserted that such varieties
are  ``built up from cones and developable varieties'' (see
[GH 79], p. 392). They gave a proof of this assertion in
the case $n=2$.
However, their assertion is not completely correct.
In a recent note [AGL],  Akivis, Goldberg,
and Landsberg   present  counter-examples to Griffiths--Harris'
assertion  when $n>2$, and in particular, they prove
that this assertion is false even for  hypersurfaces
with one-dimensional fibers.

Recently  four papers [I 98, 99a, 99b] and [IM 97] on tangentially
degenerate submanifolds (called ``developable'' in these
papers) were published. In [IM 97], the authors found the
connection between such submanifolds and solutions of
Monge-Amp\`{e}re equations, with  the foliation
of plane generators $L$ of $X$ called the Monge-Amp\`{e}re foliation.
 In [IM 97], the authors proved that the rank of the Gauss map of a
compact tangentially degenerate  $C^\infty$-hypersurface
$X \subset \mathbf{R} P^N$ is an even integer satisfying the
inequality $ \frac{r (r+3)}{2} < N, \;
r \neq 0$, and that if $r \leq 1$, then $X$ is necessarily a projective
hyperplane of $ \mathbf{R} P^N$. If $N = 3$ or $N = 5$,
then a compact tangentially degenerate  $C^\infty$-hypersurface
is a projective hyperplane.

In [I 98, 99b], Ishikawa found real algebraic cubic nonsingular
tangentially degenerate hypersurface in $\mathbf{R} P^N$ for
$N = 4, 7, 13, 25$, and in [I 99a] he studied singularities
of tangentially degenerate  $C^\infty$-hypersurfaces.

In 1997 Borisenko published the survey paper [B 97] in which he discussed
results on strongly parabolic submanifolds  and related questions in
Riemannian and pseudo-Riemannian spaces of constant curvature and in
particular, in an Euclidean space $E^N$. Among other results, he
gives a description of certain classes of submanifolds
of arbitrary codimension that are analogous to the class of
parabolic surfaces in an  Euclidean space $E^3$.
Borisenko also investigates  the
local and global metric and topological properties;
indicates conditions which imply that a submanifold of
an Euclidean space $E^N$  is cylindrical; presents
results on strongly parabolic submanifolds in pseudo-Riemannian
spaces of constant curvature, and finds the relationship
with minimal surfaces.

In the current paper we  study systematically  the differential geometry of
tangentially degenerate submanifolds of a projective space
$P^N (\mathbf{C})$. By means of the focal images, three basic types
of submanifolds are  discovered:
cones, tangentially degenerate  hypersurfaces, and torsal
submanifolds. Moreover, for  tangentially degenerate
submanifolds, a structural theorem is proven. By this theorem,
 tangentially degenerate  submanifolds that do not belong to
 the basic types are foliated into submanifolds of basic types.
 In the proof we introduce irreducible, reducible, and completely
 reducible  tangentially degenerate  submanifolds. It is found
 that cones and tangentially degenerate  hypersurfaces
are irreducible and torsal    submanifolds  are completely
reducible while all other  tangentially degenerate  submanifolds
not belonging to basic types are reducible.
Particular examples of tangentially degenerate submanifolds as well
as tangentially degenerate submanifolds of low dimensions is considered in
[AGL]. Some examples of tangentially degenerate submanifolds
can be also found in [A 87].

In this paper we apply the method of exterior forms
and moving frames of Cartan [C 45] which was often
successfully used in differential geometry.

\section{The main results}

Let $X \subset  P^N (\mathbf{C})$ be  an $n$-dimensional
 smooth submanifold  with a degenerate Gauss map
 $\gamma: X \rightarrow G (n, N), \, \gamma (x) = T_x (X), \,
 x \in X$. Suppose that  $\mbox{{\rm rank}} \,
 \gamma = r < n$. Denote by $L$ a leave of this map, $L =
 \gamma^{-1} (T_x) \subset X$.

\setcounter{theorem}{0}

\begin{theorem}
A leave $L$ of the Gauss map $\gamma$ is a subspace of $P^N, \;
\dim L = n - r = l$ or its open part.
\end{theorem}

For the proof of this theorem see [AG 93] (p. 115,  Theorem 4.1).

The foliation on $X$ with leaves $L$ is called
the {\em Monge-Amp\`{e}re foliation} (see, for example,
[D 89] or [I 98, 99b]). In this paper, we extend the leaves
of the Monge-Amp\`{e}re foliation to a projective space $P^l$
assuming that $L \sim P^l$ is a plane generator of the
submanifold $X$. As a result, we have $X = f (P^l \times M^r)$,
where $M^r$ is a  parametric variety,
and $f$ is a differentiable map $f:
P^l \times M^r \rightarrow P^N$.

However, unlike a traditional definition of the foliation
(see for example, [DFN 85], \S 29), the leaves
of the Monge-Amp\`{e}re foliation  have singularities. This is
a reason that in general its leaves are not diffeomorphic to
a standard leaf.

The tangent subspace $T_x (X)$ is fixed when a point $x$
moves along $L$. This is the reason that we  denote
it by $T_L, \, L \subset T_L$. A pair $(L, T_L)$ on $X$
 depends on $r$ parameters.

With a second-order neighborhood of a pair $(L, T_L)$,
two systems of square $r \times r$ matrices
$B^\alpha = (b_{pq}^\alpha)$ and $C_i = (c_{pi}^q), \, i = 0, 1,
\ldots, l;\, p, q = l + 1, \ldots, n; \, \alpha = n + 1, \ldots , N,$
are associated. The equation $\det (C_i x^i) = 0$ defines the set
of singular points $x = (x^i) \in L$ of the map $\gamma$.
This set is an algebraic hypersurface $F_L \subset L$ of degree $r$
which is called the {\em focus  hypersurface}.
 The equation $\det (\xi_\alpha B^\alpha) = 0$ defines the set
of singular tangent hyperplanes $\xi = (\xi_\alpha) \supset T_L$
of the map $\gamma$. This set is an algebraic hypercone $\Phi_L$ with
vertex $T_L$ which is called the {\em focus  hypercone}.

It appears that the products $H^\alpha_i = B^\alpha C_i =
(h_{ipq}^\alpha)$ are symmetric. They define on $X$ the second
fundamental forms $h_i^\alpha = h_{ipq}^\alpha \theta^p \theta^q$,
where $\theta^p$ are basis forms of the manifold $M$.
We assume that {\em for   all values of parameters
 $u = (u^p) \in M$, the system of forms $h_i^\alpha$ is regular}, i.e.,
 among them there is at least one nondegenerate
 quadratic form.

Denote by $S_L$  the osculating subspace of $X$ which is constant at
all points $x \in L$ of its generator $L$.
Its dimension is $\dim \, S_L = n + m$, where $m$ is the
number of linearly independent second
fundamental forms of $X$.

We assume that the conditions of all theorems
in this paper
are satisfied for all values of parameters
 $u \in M$.

\begin{theorem} Suppose that $l \geq 1$ and $m \geq 2$, and
  the focus hypersurfaces $F_L$ and the focus
hypercones $\Phi_L$ do not have multiple components. Then
the submanifolds $X$ is foliated into $r$ families of torses with
$l$-dimensional plane generators $L$. Each of these families
depends on $r - 1$ parameters.
\end{theorem}

For the definition of torse see Example 2 in Section 2.
A manifold $X$ described in Theorem 2 is called {\em torsal}.

\begin{theorem}
  Suppose that $l \geq 2$, and the focus
hypersurfaces $F_L$ do not have multiple components
and are indecomposable.
Then the submanifold $X$ is a hypersurface
of rank $r$ in a subspace $P^{n+1} \subset P^N$.
\end{theorem}

\begin{theorem}
  Suppose that $m \geq 2$, and the focus
hypercones $\Phi_L$ do not have multiple components
and are indecomposable.
Then the submanifold $X$ is a cone with an
 $(l-1)$-dimensional vertex and
 $l$-dimensional plane generators.
\end{theorem}

The system of matrices $B^\alpha$ and $C_i$ associated
with a submanifold $X$ is said to be {\em reducible}
if these  matrices can be simultaneously reduced to a block
diagonal form:
\begin{equation}\label{eq:1}
C_i = \mbox{{\rm diag}} \; (C_{i1}, \ldots , C_{is}), \;\;
B^\alpha = \mbox{{\rm diag}} \;(B^\alpha_{1}, \ldots , B^\alpha_{s}),
\end{equation}
where $C_{it}$ and $B^\alpha_t, \, t = 1, \ldots , s,$
are square matrices of orders $r_t$, and $r_1 + r_2 + \ldots
+ r_s = r$. If such a decomposition of matrices is not possible,
the system of matrices  $B^\alpha$ and $C_i$ is called {\em
irreducible}. If
$r_1 = r_2 = \ldots = r_s = 1$, then
the system of matrices  $B^\alpha$ and $C_i$
is called {\em completely reducible}.

A manifold $X$ with a degenerate Gauss mapping is said to be {\em reducible,
irreducible} or {\em completely  reducible} if for any values of
parameters $u \in M$ the matrices $B^\alpha$ and $C_i$  are reducible,
irreducible or completely  reducible, respectively.

\begin{theorem} Suppose that a manifold $X$  is reducible, and
its matrices  $B^\alpha_t$ and $C_{it}$ defined in $(1)$ are of
order  $r_t, \, t = 1, \ldots , s$. Then $X$  is foliated into $s$
families of  $(l + r_t)$-dimensional  submanifolds of
 rank $r_t$ with $l$-dimensional plane generators.
 For $r_t = 1$, these submanifolds are torses, and for $r_t \geq
 2$, they are irreducible submanifolds  described in Theorems
$3$ and $4$.
 \end{theorem}

\section{Examples  of submanifolds with degenerate \newline Gauss maps}

Consider a few  examples of submanifolds with a degenerate
Gauss map.

\setcounter{theorem}{0}

\examp{\label{examp:1}
For $r = 0$, a submanifold $X$ is an $n$-dimensional subspace
$P^n, \, n < N$. This submanifold is the only tangentially
degenerate  submanifold  without singularities  in $P^N$.
}

\examp{\label{examp:2} Let $Y$ be a curve of class $C^p$ in
the space $P^N$. Suppose that $P^N$ is a space
of minimal dimension containing the curve $Y$. Denote by $L_y$
the osculating subspace of order $l, \; l \leq p, \; l \leq N - 1,$
of the curve $Y$ at a point $y \in Y$. Since $\dim L_y = l$,
it follows that when a point $y$ moves along
the curve $Y$, the subspace $L_y$ sweeps  a submanifold $X$ of dimension $n
= l + 1$ in the space $P^N$. At any point $x \in L_y$, the tangent
 subspace $T_x (X)$ coincides with the  osculating subspace
 $L_y^\prime$ of order $l + 1$ of the curve $Y$,
  $\dim L_y^\prime = l + 1$, and the focus hypersurface
  in $L_y$ is the osculating subspace  ${}^\prime L_y$
  of order $l - 1$ and dimension $l - 1$  of
  the curve $Y$.   Thus $\dim X = l + 1$, and the
 manifold $X$ is tangentially degenerate of rank $r = 1$.
 Such a manifold $X$ is called a {\em torse}. Conversely, a
 submanifold of dimension $n$ and rank 1 is a torse formed
 by a family of osculating subspaces of order $n - 1$ of
 a curve of class $C^p, \, p \geq n - 1$, in the space $P^N$.

In what follows, unless otherwise stated,
we  always assume that $r > 1$.
 }

\examp{\label{examp:3}
 Suppose that $S$ is a subspace of the space $P^N, \;
\dim S = l - 1$, and $T$ is its complementary subspace,
$\dim \, T = N -l, \;  T \cap S = \emptyset$.
Let $Y$ be a smooth  tangentially nondegenerate
submanifold of the subspace $T$, $\dim Y = \mbox{{\rm rank}} \; Y
 = r < N -l$. Consider an $r$-parameter
family of $l$-dimensional subspaces $L_y = S \wedge y, \; y \in Y$. This
manifold is a cone $X$ with vertex $S$ and the director manifold
$Y$. The subspace $T_x (X)$ tangent to the cone $X$ at a point $x
\in L_y$ is defined by its vertex $S$ and the subspace $T_y (Y),
\; T_x (X) = S \wedge T_y (X)$, and  $T_x (X)$ remains fixed when a point $x$
moves in the subspace $L_y$. As a result, the cone $X$ is a
tangentially degenerate submanifold of dimension $n = l + r$ and
 rank $r$, with plane generators $L_y$ of dimension $l$.
}

\examp{\label{examp:4} Let $T$ be a subspace of dimension $n + 1$
 in $P^N$, and let $Y$ be an $r$-parameter family of hyperplanes
 $\xi$ in general position in $T, \, r < n$. Such a family has an
 $n$-dimensional envelope $X$ that is a tangentially degenerate
 submanifold  of dimension $n$ and rank $r$ in the subspace $T$.
 It is foliated into an $r$-parameter family of plane generators
 $L$  of dimension $l = n -r$ along which the tangent subspace
 $T_x (X), \, x \in L,$ is fixed and coincides with a hyperplane
 $\xi$ of the family in question. Thus $X$ is a  tangentially degenerate
 hypersurface  of rank $r$ with  $(n-r)$-dimensional flat
 generators $L$ in the space $T$.
}

\examp{\label{examp:5}
If $\mbox{{\rm rank}}\, X = \dim \, X = n$, then $X$ is a submanifold
of complete rank. $X$ is also called a
{\em tangentially nondegenerate} in the space $P^N$.
}

\section{Application of the duality principle}

By the duality principle, to a point $x$ of a projective space $P^N$,
there corresponds a hyperplane $\xi$. A set of hyperplanes
of space $P^N$ forms the dual  projective space $(P^N)^*$
of the same dimension $N$. Under this correspondence,
to a subspace $P$ of dimension $p$,
there corresponds a subspace $P^* \subset (P^N)^*$ of
dimension $N - p - 1$. Under the dual map, the incidence
of subspaces is preserved, that is, if $P_1 \subset P_2$, then
$P_1^* \supset P_2^*$.

In the space $P^N$, if a point $x$ describes a  tangentially
nondegenerate manifold
$X, \; \dim X = r$,  then in general, the hyperplane $\xi$
corresponding to $x$ envelopes a hypersurface $X^*$ of rank $r$
with $(N - r - 1)$-dimensional generators. The dual map sends
a smooth manifold $X \subset P^N$ of dimension $n$ and rank $r$
with plane generators $L$ of dimension $l = n - r$ to
a manifold $X^*$  of dimension $n^* = N - l - 1$ and the same rank $r$
with plane generators $L^*$ of dimension $l^* = N - n - 1$.
Under this map, to a tangent subspace $T_x (X)$ of the submanifold $X$
there corresponds the plane generator $L^*$, and to a plane
generator $L$ there corresponds the tangent subspace of
the  manifold $X^*$.

Note that Ein [E 85, 86] applied the duality
principle for studying tangentially degenerate varieties with
small dual varieties.

Let us determine which manifolds correspond to
tangentially degenerate manifolds considered in Examples 2--4.
To a cone $X$ of rank $r$ with vertex $S$ of dimension $l - 1$
(see Example 3), there corresponds a  manifold
$X^*$ lying in the subspace $T = S^*, \; \dim T = N - l$.
Since $\dim X^* = n^* = N - l - 1$, the  manifold
$X^*$ is a hypersurface of rank $r$ in the subspace $T$.
Such a hypersurface was considered in Example 4.
It follows that Examples 3 and 4 are mutually dual
one to another.

Suppose that the subspace $T$ containing a hypersurface $X$ of
rank $r$ coincides with the space $P^N$. Then
under the dual map in $T$, to the   hypersurface $X$
there corresponds an $r$-dimensional tangentially nondegenerate
manifold $X^*$.

Under the dual map, in $P^N$ to torses (see Example 2) there
correspond torses enveloping a one-parameter family of hyperplanes
$\eta$---the images of points of the curve $Y$. These torses
are of dimension $N - 1$ and rank 1, and their plane generators
are of dimension $l = N - 2$.

Note that the dual map allows us to construct tangentially degenerate
manifolds in the space $P^{n+1}$ from the general
tangentially nondegenerate manifolds.

\examp{\label{examp:6}
 Consider a set of conics in a projective plane $P^2$. They are defined
by the equation
\begin{equation}\label{eq:2}
a_{ij} x^i x^j = 0, \;\; i, j = 0, 1, 2,
\end{equation}
where the coefficients $a_{ij}$ are defined up to
a constant factor, and $a_{ij} =  a_{ji}$.
Thus to the curve (2) there corresponds a point
of a five-dimensional projective space $P^{5*}$
whose coordinates coincide with the coefficients of
equation (2).

Among the conics defined by equation (2) there are singular curves  which
decompose into two straight lines. Such conics are distinguished
by the condition
\begin{equation}\label{eq:3}
\det \, (a_{ij}) = 0.
\end{equation}
Equation (3) is of degree three and defines
a hypersurface  in $P^{5*}$ which
is called the {\em cubic symmetroid.}

We prove now that {\em the hypersurface defined by equation
$(3)$
is tangentially degenerate of rank two  and bears
a two-parameter family of two-dimensional plane generators.}
To this end, we write equation (3) in the form
\begin{equation}\label{eq:4}
F = \det \pmatrix{a_{00} & a_{01} & a_{02} \cr
                  a_{10} & a_{11} & a_{12} \cr
                  a_{20} & a_{21} & a_{22} \cr} = 0.
\end{equation}
The equation of tangent subspaces of the hypersurface (4) has the
form
$$
\xi^{ij} a_{ij} = 0,
$$
where $\xi^{ij} = \displaystyle \frac{\partial F}{\partial a_{ij}}, \; \xi^{ij} =
\xi^{ji},$ are the cofactors of the entries $a_{ij}$ in the
determinant in equation (4). The quantities $\xi^{ij}$ are coordinates
in the projective space $P^5$, which is dual to the space
$P^{5*}$  where we defined the cubic symmetroid (4).

Consider the symmetric matrix $\xi = (\xi^{ij})$. A
straightforward computation shows that the rank of this matrix
is equal to one. Thus the entries of the matrix $\xi$ can be
represented in the form
\begin{equation}\label{eq:5}
\xi^{ij} = \xi^{i}\xi^{j},
\end{equation}
where
$\xi^{i}$ are projective coordinates of a point in the
projective plane $P^2$. But equations (5)
 define the Veronese surface
 in the space $P^5$. Since the Veronese surface  is
 of dimension two, its points depend on two
 affine parameters $u = \frac{\xi^1}{\xi^0}$ and
 $v = \frac{\xi^2}{\xi^0}$. As a result, the tangent
 hyperplanes of the cubic symmetroid (4) depend on
 two parameters, and this hypersurface is of rank two.

To the plane generators  of the hypersurface (4),
there correspond two-parameter families
of conics in $P^2$ decomposing into
pairs of straight lines with a common
intersection point.

Thus, {\em the cubic symmetroid $(4)$ and the Veronese surface  $(5)$
are mutually dual submanifolds of a five-dimensional projective
space.}
}

\section{Basic equations and focal images}

We  study  tangentially degenerate  submanifolds applying
the method of moving frames in a projective space $P^N$.
In $P^N$, we consider a manifold of projective frames
$\{A_0, A_1, \ldots , A_N\}$. On this manifold
\begin{equation}\label{eq:6}
 d A_u =  \omega_u^v A_v, \;\; u, v, = 0, 1, \ldots , N,
\end{equation}
where the sum $\omega_u^u = 0$. The 1-forms
$\omega_u^v$ are linearly expressed in terms of
the differentials of parameters of the group
of projective transformations of the space $P^N$. The total
number of these parameters is $N^2 - 1$. These 1-forms
satisfy the structure equations
\begin{equation}\label{eq:7}
 d  \omega_u^v  =  \omega_u^w \wedge  \omega_w^v
\end{equation}
of the space $P^N$ (see, for example, [AG 93], p. 19).
Equations (7) are the conditions of complete integrability
of equations (6).

Consider a  tangentially degenerate submanifold $X \subset P^N,
\; \dim X = n, \linebreak \mbox{{\rm rank}}\; X = r \leq n$. In addition,
as  above, let $L$ be a rectilinear generator of the manifold
$X, \; \dim L = l$; let $T_L, \; \dim T_L = n,$ be the tangent
subspace to $X$ along the generator $L$, and let $M$ be a base
manifold for $X, \; \dim M = r$. Denote by $\theta^p, \; p = l + 1,
\ldots, n$,  basis forms on the variety $M$. These forms satisfy
the structure equations
\begin{equation}\label{eq:8}
 d  \theta^p  =  \theta^q \wedge  \theta_q^p, \;\; p, q = l + 1,
 \ldots , n,
\end{equation}
of the variety $M$. Here $\theta_q^p$ are 1-forms defining
transformations of first-order frames on $M$.

For a point $x \in L$, we have $d x \in T_L$. With $X$, we
associate a bundle of projective frames $\{A_i, A_p, A_\alpha\}$
such that $A_i \in L, \; i = 0, 1, \ldots, l; \; A_p \in T_L, \;
p = l + 1, \ldots , n$. Then
\begin{equation}\label{eq:9}
\renewcommand{\arraystretch}{1.3}
\left\{
\begin{array}{ll}
dA_i =  \omega_i^j A_j + \omega_i^p A_p, \\
dA_p =  \omega_p^i A_i + \omega_p^q A_q +
 \omega_p^\alpha A_\alpha, \; \alpha = n + 1, \ldots, N.
\end{array}
\right.
\renewcommand{\arraystretch}{1}
\end{equation}
It follows from the first equation of (9) that
\begin{equation}\label{eq:10}
  \omega_i^\alpha = 0.
\end{equation}
Since for $\theta^p = 0$ the subspaces $L$ and $T_L$ are fixed, we
have
\begin{equation}\label{eq:11}
  \omega_i^p = c^p_{qi} \theta^q, \;\;   \omega_p^\alpha
= b^\alpha_{qp} \theta^q.
\end{equation}

Since the manifold of leaves $L \subset X$ depends on $r$
essential parameters, the rank of the system of 1-forms
$\omega^p_i$ is equal to $r$, $\mbox{{\rm rank}} \; (\omega^p_i) =
r$. Similarly, we have $\mbox{{\rm rank}} \; (\omega^\alpha_p) =
r$.

 Denote by
$C_i$ and $B^\alpha$ the  $r \times r$ matrices occurring in
equations (11):
$$
C_i = (c^p_{qi}), \;\; B^\alpha = (b_{qp}^\alpha).
$$
These matrices are defined in a second-order neighborhood
of the submanifold $X$.

Exterior differentiation of equations (10) by
means of  structure equations (7) leads to the exterior
quadratic equations
$$
\omega_p^\alpha \wedge \omega_i^p = 0.
$$
Substituting expansions (11) into the last equations, we find that
\begin{equation}\label{eq:12}
b^\alpha_{qs} c^s_{pi} = b^\alpha_{ps} c^s_{qi}.
\end{equation}
Equations (10), (11), and (12) are called the {\em basic equations}
in the theory of tangentially degenerate submanifolds.

Relations (12) can be written in the matrix form $$ (B^\alpha
C_i)^T = (B^\alpha C_i), $$ i.e., the matrices $$ H^\alpha_i =
B^\alpha C_i = (b^\alpha_{qs} c^s_{pi}) $$ are symmetric.

Let $x = x^i A_i$ be an arbitrary point of a leaf $L$. For such a
point we have $$ d x = (dx^i + x^j \omega_j^i) A_i + x^i
\omega_i^p A_p. $$ It follows that $$ d x \equiv (A_p c^p_{qi}
x^i) \theta^q \pmod{L}. $$ The tangent subspace $T_x$ to the
manifold $X$ at a point $x$ is defined by the points $A_i$ and $$
\widetilde{A}_q (x) = A_p c^p_{qi} x^i, $$ and therefore $T_x
\subset T_L$.

A point $x$ is a {\em regular} point of a leaf $L$ if
$T_x = T_L$. Regular points are determined by the condition
\begin{equation}\label{eq:13}
J (x) = \det (c^q_{pi} x^i) \neq 0.
\end{equation}
If $J (x) = 0$ at a point $x$, then $T_x$ is a proper subspace
of $T_L$, and a point $x$ is said to be a {\em singular point}
of a leaf $L$.

The determinant (13) is the Jacobian of the map $f: P^l \times
M^r \rightarrow P^N$.
Singular points of a leaf $L$ are determined by the
condition $J (x) = 0$. In a leaf $L$, they form an algebraic
submanifold of dimension $l - 1$ and degree $r$. This
hypersurface (in $L$) is called the {\em focus hypersurface}
and is denoted by $F_L$.
By (13), the equations $J (x) = 0$ of the focus hypersurface
on the plane generator $L$ of the manifold $X$ can be written as
\begin{equation}\label{eq:14}
 \det (c^q_{pi} x^i) = 0.
\end{equation}

We calculate now the second differential of
a point $x \in L$:
$$
d^2 x \equiv A_\alpha \omega_s^\alpha \omega_i^s x^i \pmod{T_x}.
$$
This expression is the {\em second fundamental quadratic form}
of the manifold $X$:
\begin{equation}\label{eq:15}
II_x = A_\alpha \omega_s^\alpha \omega_i^s x^i  =
A_\alpha b^\alpha_{ps} c^s_{qi} x^i  \theta^p \theta^q.
\end{equation}

Suppose that $\xi = \xi_\alpha x^\alpha = 0$ is the tangent
hyperplane to $X$ at $x \in L, \; \xi \supset T_L$. Then
$$
(\xi, II_x) = h_{pq} (\xi, x)\theta^p \theta^q,
$$
where
$$ h_{pq} (\xi,
x) = \xi_\alpha b^\alpha_{ps} c^s_{qi} x^i, \;\;\;\; h_{pq} = h_{qp},
$$
 is the {\em second fundamental quadratic form} of the manifold
$X$ at $x$ with respect to the hyperplane $\xi$.
Since at regular points $x \in L$ condition (13) holds, the rank of
this matrix is the same as the rank of the matrix
\begin{equation}\label{eq:16}
B (\xi) = (\xi_\alpha b^\alpha_{pq}) = \xi_\alpha B^\alpha,
\end{equation}
and this rank is the same at all  regular points $x \in L$.

We call a tangent hyperplane $\xi$ {\em singular} if
\begin{equation}\label{eq:17}
\det (\xi_\alpha b^\alpha_{pq}) = 0,
\end{equation}
i.e., if the rank of the matrix (16) is reduced. Condition (17) is
an equation of degree $r$ with respect to the tangential
coordinates $\xi_\alpha$ of the hyperplane $\xi$. This condition
defines an algebraic hypercone whose vertex is the tangent
subspace $T_L$. This hypercone is called the {\em focus hypercone}
and denoted by $\Phi_L$ (see [AG 93], p. 119).

The determinant $\det (\xi_\alpha b^\alpha_{pq})$ in the
left-hand side of equation (17) is the Jacobian of the dual
map $f^*: L^* \times M^r \rightarrow (P^N)^*, \;
f^* (L^* \times M^r) = X^* \subset  (P^N)^*$, where $X^*$ is
a manifold dual to $X$ and $L^*$ is a bundle of hyperplanes of
the space $P^N$ passing through the tangent subspace $T_L$
of the manifold $X, \; \dim \, L^* = N - n - 1$.

The focus hypersurface $F_L \subset L$ and
the focus hypercone $\Phi_L$ with vertex  $T_L$ are called
{\em focal images} of the manifold $X$ with a degenerate Gauss
map.

Note that under the passage from the manifold $X \subset P^N$
to its dual manifold $X^* \subset (P^N)^*$, the systems of
square matrices $C_i$ and $B^\alpha$ as well  as
the focus hypersurfaces $F_L$ and the focus cones $\Phi_L$
exchange their roles.

Since
$$
d^2 x \equiv A_\alpha b_{qs}^\alpha c_{pi}^s x^i
\theta^p \theta^q \pmod{T_L, x \in L},
$$
the points
\begin{equation}\label{eq:18}
A_{pq} = A_\alpha b^\alpha_{qs} c_{pi}^s x^i, \;\; A_{pq} =
A_{qp},
\end{equation}
together with the points $A_i$ and $A_p$ define the
osculating subspace $T^2_L (X)$. Its dimension is
$$
\dim T_L^2 (X) = n + m,
$$
where $m$ is the number of linearly independent points
among the points
$A_{pq}, \; m \leq \frac{r (r + 1)}{2}$. But since
at a regular point $x \in L$  condition (13) holds,
the number $m$ is the number of linearly independent points
among the points
$$
\widetilde{A}_{pq} = A_\alpha b^\alpha_{pq}.
$$
We  also use the notation $S_L$ for
 the osculating space $T_L^2 (X)$.

On a generator $L$ of the manifold $X$, consider the
system of equations
\begin{equation}\label{eq:19}
 c^q_{pi} x^i = 0.
\end{equation}
Its matrix $C = (c^q_{pi})$ has $r^2$ rows and $l + 1$ columns.
Denote the rank of this matrix by $m^*$. If $m^* < l + 1$,
then  system (19) defines a subspace $K_L$ of dimension $k = l - m^*$
in $L$. This subspace belongs to the focus hypersurface
$F_L$ defined by equation (14). If $l > m^*$, then
the hypersurface $F_L$ becomes a cone with vertex $K_L$. We  call
the subspace $K_L$ the {\it characteristic subspace} of the
generator $L$.

Note also that by the duality principle in $P^N$, the osculating
subspace $S_L$ and the characteristic subspace $K_L$ constructed for
a pair $(L, T_L)$ correspond one to another.

\section{Proof of Theorem 2}

\setcounter{theorem}{5}

\begin{lemma} Suppose that $l \geq 1$, and the focus
hypersurface $F_L \subset L$ does not have multiple components. Then all
matrices $B^\alpha$ can be simultaneously diagonalized,
$B^\alpha = (\mbox{{\rm diag}} \;b_{pp}^\alpha)$, and the focus
hypercone $\Phi_L$ decomposes into $r$ bundles  of hyperplanes $\Phi_p$
in $P^N$ whose axes are $(n+1)$-planes $T_L \wedge  B_p$, where
$B_p = b_{pp}^\alpha A_\alpha$ are points located outside of the tangent
subspace $T_L$. The dimension $n + m$ of the osculating subspace $S_L$ of the
manifold $X$ along a generator $L$ does not exceed $n + r$.
\end{lemma}

{\sf Proof.}  Since the hypersurface $F_L \subset L$ does not have multiple
components, a general  straight line $\lambda$
lying in $L$ intersects $F_L$ at $p$
distinct points. We place the vertices $A_0$ and $A_1$ of our
moving frame onto the line $\lambda$. By (14), the coordinates of
the common points of $\lambda$ and $F_L$ are defined by the equation
$$
\det (c_{p0}^q x^0 + c_{p1}^q x^1 ) = 0.
$$
Suppose that $A_0 \notin F_L$. Then $\det \;(c_{p0}^q) \neq 0$, and it
is easy to prove that the matrices $C_0$ and $C_1$ can be
simultaneously  diagonalized,  $C_0 = (\delta_q^p)$ and  $C_1 = (\mbox{{\rm diag}}
\;\; c_{p1}^p)$. Since the common points of $\lambda$ and $F$
are not multiple points, we have $c_{p1}^p \neq c_{q1}^q$ for $p \neq q$.

Next we write equations (12) for $i = 0, 1$:
$$
b^\alpha_{pq} = b^\alpha_{qp}, \;\; b^\alpha_{qp} c_{p1}^p = b^\alpha_{pq}
c_{q1}^q.
$$
Since  $c_{p1}^p \neq c_{q1}^q$ , it follows that $b^\alpha_{pq}
= 0$ for $p \neq q$. As a result, all  matrices $B^\alpha$
can be simultaneously diagonalized,
$B^\alpha = (\mbox{{\rm diag}} \;b_{pp}^\alpha)$. Equation (17)
takes the form
$$
\prod_p (\xi_\alpha b^\alpha_{pp}) = 0,
$$
and the focus hypercone $\Phi_L$
decomposes into $r$ bundles of hyperplanes  $\Phi_p$ in $P^N$ whose axes are
$(n+1)$-planes $T_L \wedge  B_p$, where $B_p = b_{pp}^\alpha A_\alpha$
are points located outside of the tangent subspace $T_L$.
 The osculating subspace $S_L$ of the
manifold $X$ along a generator $L$  is the span
of the tangent subspace $T_L$ and the points $B_{l+1}, \ldots,
B_n$. Thus, the dimension of this subspace does not exceed $n + r$.
 \rule{3mm}{3mm}

\begin{lemma} Suppose that $l \geq 1, \; m \geq 2$,
the focus hypersurfaces $F_L \subset L$ do not have
multiple components, and all the bundles
$\Phi_p$ into which the hypercone $\Phi_L$ decomposes
are of multiplicity one. Then the focus hypersurface $F_L$
decomposes into $r$ hyperplanes $F_p \subset L$.
\end{lemma}

{\sf Proof.} Consider the matrix
\begin{equation}\label{eq:20}
B = (b^\alpha_{pp})
\end{equation}
composed from the eigenvalues of the matrices $B^\alpha$.
Matrix (20) has $r$ columns and $m$ independent rows, $m \leq r$.
Write equations (12) for the diagonal matrices $B^\alpha$:
\begin{equation}\label{eq:21}
b^\alpha_{qq} c^q_{pi} = b^\alpha_{pp} c^p_{qi}.
\end{equation}

Since the matrix $B$ has $m$ linearly independent
columns and $m \geq 2$,  it follows from equation (21)
that $c_{pi}^q = 0$ for $p \neq q$, and all the
matrices $C_i$ and $B^\alpha$ can be
simultaneously diagonalized. Now the equation of the
focus hypersurface $F_L$ takes the form
$$
\prod_{p = l + 1}^n c_{pi}^p x^i = 0,
$$
and the  hypersurface $F_L$ decomposes into $r$ hyperplanes $F_p$
 defined in $L$ by the equation $ c_{pi}^p x^i =
0$. Since by the conditions of Theorem  2,
 the focus hypersurface $F_L$ does not have multiple components,
all  hyperplanes $F_p$ are distinct.  \rule{3mm}{3mm}

Consider a rectangular $r \times (l + 1)$ matrix
\begin{equation}\label{eq:22}
C = (c^p_{pi})
\end{equation}
formed by the eigenvalues of the matrix $C_i$.

{\sf Proof of Theorem 2}.
From the conditions of Theorem 2 and Lemmas 6 and 7, it follows
that the matrices $C_i$ and $B^\alpha$ can be simultaneously
diagonalized:
$$
 C_i = \mbox{{\rm diag}} \; (c^{l+1}_{l+1, i}, \ldots , c^n_{ni}), \;\;
 B^\alpha = \mbox{{\rm diag}} \; (b^\alpha_{l+1, l+1}, \ldots , b^\alpha_{nn}).
$$

This implies that formulas (11) take the form:
\begin{equation}\label{eq:23}
  \omega_i^p = c^p_{pi} \theta^p, \;\;   \omega_p^\alpha
= b^\alpha_{pp} \theta^p,
\end{equation}
where $\theta^p$ are basis forms on the  parametric variety $M$,
 and there is no summation over the index $p$.
Exterior differentiation of equations (23) leads to the following
exterior quadratic equations:
\begin{equation}\label{eq:24}
  \nabla c^p_{pi} \wedge \theta^p + c^p_{pi}d \theta^p
  + \sum_{q \neq p}  c^q_{qi} \omega_q^p  \wedge \theta^q = 0,
\end{equation}
\begin{equation}\label{eq:25}
  \nabla b^\alpha_{pp} \wedge \theta^p + b^\alpha_{pp} d\theta^p
  - \sum_{q \neq p}  b^\alpha_{qq} \omega_p^q  \wedge \theta^q = 0,
\end{equation}
where
$$
 \nabla c^p_{pi} =  d c^p_{pi} -  c^p_{pj} \omega_i^j +  c^p_{pi}
 \omega_p^p, \;\;
 \nabla b^\alpha_{pp} =  d  b^\alpha_{pp}
 + b^\beta_{pp} \omega_\beta^\alpha - b^\alpha_{pp}  \omega_p^p.
$$

Suppose now that at least one of the matrices $B$ and $C$
defined by equations (20) and (22) has mutually linearly independent
columns. Then it follows from equations (24) and (25) that
$$
\renewcommand{\arraystretch}{1.3}
\left\{
\begin{array}{ll}
d\theta^p \equiv 0 \pmod{\theta^p}, \\
\omega_q^p \wedge \theta^q  \equiv 0 \pmod{\theta^p}, \\
\omega_p^q \wedge \theta^q  \equiv 0 \pmod{\theta^p}.
\end{array}
\right.
\renewcommand{\arraystretch}{1}
$$
This implies that
\begin{equation}\label{eq:26}
d\theta^p  = \theta^p \wedge \theta^p_p, \;\;
\omega_p^q  = l_{pp}^q \theta^p +  l_{pq}^q \theta^q.
\end{equation}
We write the expression of $dA_p$ from (9)
in more detail:
\begin{equation}\label{eq:27}
dA_p =  \omega_p^i A_i + \omega_p^p A_p + \sum_{q \neq p}
\omega_p^q A_q + \omega^\alpha_p A_\alpha.
\end{equation}

Let the index $p$ be fixed and $q \neq p; \,
p, q = l + 1, \ldots , n$. Consider a submanifold
$X_p \subset X$ defined by the equations
\begin{equation}\label{eq:28}
\theta^q  = 0, \;\; q \neq p.
\end{equation}
The first equation of (26) implies that $\dim X_p = 1$.
By (28), (9) and (27), we find that
\begin{equation}\label{eq:29}
dA_i =  \omega_i^j A_j + c_{pi}^p \theta^p A_p,
\end{equation}
\begin{equation}\label{eq:30}
dA_p =  \omega_p^i A_i + \omega_p^p A_p +
B_p \theta^p,
\end{equation}
where $B_p = \sum_{q\neq p} l^q_{pp} A_q +  b^\alpha_{pp}  A_\alpha$.
Since the points $A_i$ are basis points of a generator
$L$ of the manifold $X$, equations (29) and (30) prove that the
subspace $L \wedge A_p$ is  tangent to the submanifold
$X_p$ at all regular points of the  generator
$L$, and the subspace  $L \wedge A_p \wedge B_p$
is the osculating subspace
of $X_p$. Singular points of a generator $L$  of the manifold $X_p$
are determined by the equations $c_{pi}^p x^i = 0$ and
form a hyperplane in $L$ (see Example  2 in Section  2).
Thus,  the submanifold $X_p$
is a torse with $l$-dimensional plane generators. Hence
 the manifold $X$ is foliated into $r$ families of torses,
 each of these families depends on $r - 1$ parameters $u^q$,
 and the forms $\theta^q$ are expressed in terms of the
 differentials of these parameters, $\theta^q = l^q d u^q$.
\rule{3mm}{3mm}

\section{Proof of Theorems 3 and 4}

\begin{lemma} Suppose that $l \geq 2$,  the
focus hypersurfaces $F_L \subset L$ do not have
multiple components and are indecomposable.
Then the hypercone $\Phi_L$
is an $r$-multiple bundle of hyperplanes with
an $(n+1)$-dimensional vertex in $P^N$.
\end{lemma}

{\sf Proof.} From the conditions of the lemma and equation (21)
it follows that the columns of the matrix $B$ are linearly dependent,
and thus all points $B_p = b_{pp}^\alpha A_\alpha$ defined by
the columns of this matrix lie in  a subspace of dimension $n + 1$
defined by the tangent subspace $T_L$ and one of the points $B_p$.

Therefore all bundles of hyperplanes $\Phi_p$ into which the
hypercone $\Phi_L$ decomposes coincide. \rule{3mm}{3mm}

We  formulate the lemma which is dual to Lemma 8.

\begin{lemma} Suppose that $m \geq 2$ and the focus
hypercones $\Phi_L$ with their vertices $T_L$ do
not have multiple components and are indecomposable.
Then the focus hypersurface $F_L \subset L$
is an $r$-multiple   hyperplane in $L$.
\end{lemma}

{\sf Proof of Theorem 3}. From the conditions of
Theorem 3 and Lemma 8, it follows
that all pairs of columns of the matrix $B$ are linearly dependent. Thus
all matrices $B^\alpha$ associated with $X$ are mutually linearly
dependent. Hence, we have
\begin{equation}\label{eq:31}
b_{pq}^\alpha = b^\alpha b_{pq}.
\end{equation}
Now conditions (12) take the form
\begin{equation}\label{eq:32}
b_{ps} c^s_{qi} = b_{qs} c^s_{pi}.
\end{equation}
Although the matrix $B = (b_{pq})$ is still can be diagonalized,
in general, the matrices $C_i$ do not possess this property.
Thus we write equations (11) in the form
\begin{equation}\label{eq:33}
  \omega_i^p = c^p_{qi} \theta^q, \;\;   \omega_p^\alpha
= b^\alpha b_{pq} \theta^q.
\end{equation}

Exterior differentiation of equations (33) and applying  structure
equations (7) of  the  parametric variety $M$, we obtain the following
exterior quadratic equations:
\begin{equation}\label{eq:34}
  \nabla c^p_{qi} \wedge \theta^q = 0, \;\;
    (b_{pq} \nabla b^\alpha + b^\alpha \nabla b_{pq}) \wedge \theta^q = 0,
\end{equation}
where
$$
\renewcommand{\arraystretch}{1.3}
\left\{
\begin{array}{ll}
 \nabla c^p_{qi} =  d c^p_{qi} -  c^p_{qj} \omega_i^j +  c^s_{qi}
 \omega_s^p - c_{si}^p \theta_p^s, \\
 \nabla b^\alpha =  d  b^\alpha + b^\beta \omega_\beta^\alpha, \\
  \nabla b_{pq} =  d  b_{pq} -  b_{sq} \omega_p^s - b_{ps}  \theta_q^s.
\end{array}
\right.
\renewcommand{\arraystretch}{1}
$$

As we  noted already in the proof of Theorem 2, if
the point $A_0$ does not belong to the focus
hypersurface $F_L$ of a generator $L$, then the
matrix $(c^p_{q0})$ is nonsingular, and by
means of a frame transformation in the tangent space
$T_u (M)$ to the variety $M$ we can reduce this matrix
to the form $(c^p_{q0}) = (\delta_q^p)$.
As a result, equation (32) corresponding to
$i = 0$ takes the form
\begin{equation}\label{eq:35}
b_{pq}  = b_{qp}.
\end{equation}
Hence the matrix $B = (b_{pq})$ becomes symmetric.
This matrix is the matrix of the second fundamental
form $II$ of the manifold $X$ at the point $A_0$.
Since $A_0$ is a regular point of $X$, the matrix
$B$ is nonsingular, $\det (b_{pq}) \neq 0$.

Next, we find from the second equation of (34) that
$$
   \nabla b^\alpha = b_p^\alpha \theta^p, \;\;
    \nabla b_{pq} = b_{pqs} \theta^s.
$$
Substituting these expansions into
equation (34) and equating to 0 the coefficients
in independent products $\theta^s \wedge \theta^q$,
we find that
$$
  b_{pq}   b_s^\alpha -  b_{ps}   b_q^\alpha
  +  b^\alpha (b_{pqs} - b_{psq}) = 0.
$$
Contracting these equations with the matrix
$B^{-1} = (b^{pq})$ which is the inverse of $B$,
we find that
$$
(r-1) b_s^\alpha +   b^\alpha (b_{pqs} - b_{psq}) b^{pq} = 0.
$$
Since by the theorem hypothesis $r \neq 1$, it follows that
\begin{equation}\label{eq:36}
b_s^\alpha  = b^\alpha b_s,
\end{equation}
where
$$
b_s = \frac{1}{r-1} (b_{pqs} - b_{psq}) b^{pq}.
$$

Next consider equations (17) of the focus hypercone
$\Phi_L$ of the manifold $X$. By (31), this equation
becomes
$$
(b^\alpha \xi_\alpha)^r \det (b_{pq}) = 0,
$$
and since $ \det (b_{pq}) \neq 0$, the
 hypercone $\Phi_L$ becomes an $r$-multiple bundle of hyperplanes.
 The axis of this bundle is a subspace of dimension
 $n + 1$ which is the span of the tangent subspace $T$ and the
 point $b^\alpha A_\alpha$.

 Let us prove that if the parameters $u$ on the variety $M$
 vary, this subspace is fixed. In fact, the basis points of
the subspace $T$ are the points $A_i$ and $A_p$, and by
(9) and (36), we have
$$
\renewcommand{\arraystretch}{1.3}
\left\{
\begin{array}{ll}
dA_i =  \omega_i^j A_j + \omega_i^p A_p, \\
dA_p =  \omega_p^i A_i + \omega_p^q A_q +
b_{pq} \theta^q \cdot b^\alpha A_\alpha.
\end{array}
\right.
\renewcommand{\arraystretch}{1}
$$
If we differentiate the point $b^\alpha A_\alpha$ and
apply equation (36), we find that
$$
d (b^\alpha A_\alpha) \equiv b_s \theta^s \cdot b^\alpha A_\alpha
\pmod{A_i, A_p}.
$$
This proves our last assertion.

Thus the subspace $P^{n+1} = T \wedge b^\alpha A_\alpha$
is fixed when the tangent subspace $T$ moves along $X$, and
$X$ is a hypersurface in  the subspace $P^{n+1}$.
 \rule{3mm}{3mm}

{\sf Proof of Theorem 4}. Theorem 4 is dual to Theorem 3 and can
be proved by applying Lemma 9 in the same way as we used
Lemma 8 to prove Theorem 3.  \rule{3mm}{3mm}

\section{Proof of Theorem 5}

Equations (14) and (17) of the focal images imply the following
lemma.

\begin{lemma} If  the matrices $C_i$ and $B^\alpha$  of  a manifold $X$
can be reduced to the form $(1)$, then each of its focus hypersurfaces
$F_L \subset L$ decomposes into $s$ components $F_t$ of
dimension $l - 1$ and degree $r_1, r_2, \ldots , r_s$,
and  each of its focus hypercones $\Phi_L$ decomposes into $s$ hypercones
$\Phi_t$ of the same degrees $r_1, r_2, \ldots , r_s; \, r_1 + r_2 + \ldots
+ r_s = r,$ and with the same vertex $T$. In particular, if
$r_1 = r_2 = \ldots = r_s = 1$, then
a  focus hypersurface $F_L$ decomposes  into $r$ hyperplanes,
and  a focus hypercone $\Phi_L$ decomposes into $r$ bundles
of hyperplanes with $(n+1)$-dimensional axes.
\end{lemma}

{\sf Proof of Theorem 5}.

 We  prove Theorem 5 assuming that the index $t$
takes only two values, $t = 1, 2, \; r = r_1 + r_2$, and the
indices $p$ and $q$ have the following values:
$$
p_1, q_1 = l + 1, \ldots , l + r_1, \;\;
p_2, q_2 = l + r_1 + 1, \ldots , n.
$$
Then equations (11) become
\begin{equation}\label{eq:37}
\renewcommand{\arraystretch}{1.3}
\left\{
\begin{array}{ll}
  \omega_i^{p_1} = c^{p_1}_{q_1i} \theta^{q_1}, &
  \omega_{p_1}^\alpha = b^\alpha_{p_1 q_1} \theta^{q_1}, \\
  \omega_i^{p_2} = c^{p_2}_{q_2i} \theta^{q_2}, &
  \omega_{p_2}^\alpha = b^\alpha_{p_2 q_2} \theta^{q_2}.
\end{array}
\right.
\renewcommand{\arraystretch}{1}
\end{equation}

Exterior differentiation of equations (37) gives
\begin{equation}\label{eq:38}
  \nabla c^{p_1}_{q_1i} \wedge \theta^{q_1} +
(c^{s_2}_{q_2i}   \omega^{p_1}_{s_2} - c^{p_1}_{s_1 i} \theta_{q_2}^{s_1})
 \wedge  \theta^{q_2} = 0,
\end{equation}
 \begin{equation}\label{eq:39}
  \nabla b^\alpha_{p_1 q_1} \wedge \theta^{q_1} -
(b_{s_2 q_2}^\alpha   \omega_{p_1}^{s_2} + b_{p_1 s_1}^\alpha
\theta_{q_2}^{s_1} \wedge \theta^{q_2} = 0,
\end{equation}
\begin{equation}\label{eq:40}
  \nabla c^{p_2}_{q_2i} \wedge \theta^{q_2} +
(c^{s_1}_{q_1i}   \omega^{p_2}_{s_1} - c^{p_2}_{s_2 i} \theta_{q_1}^{s_2})  \wedge  \theta^{q_1} = 0,
\end{equation}
\begin{equation}\label{eq:41}
  \nabla b^\alpha_{p_2 q_2} \wedge \theta^{q_2} -
(b_{s_1 q_1}^\alpha   \omega_{p_2}^{s_1} + b_{p_2 s_2}^\alpha
\theta_{q_1}^{s_2} \wedge \theta^{q_1} = 0,
\end{equation}
where
$$
\renewcommand{\arraystretch}{1.3}
\left\{
\begin{array}{ll}
\nabla c^{p_1}_{q_1i} = d c^{p_1}_{q_1i} - c^{p_1}_{q_1j}
\omega_i^j + c^{s_1}_{q_1i} \omega_{s_1}^{p_1}
- c^{p_1}_{s_1i} \theta^{s_1}_{q_1}, \\
\nabla b^\alpha_{p_1 q_1} = d b^\alpha_{p_1 q_1} + b^\beta_{p_1 q_1}
\omega_\beta^\alpha - b^\alpha_{s_1 q_1} \omega^{s_1}_{p_1}
- b^\alpha_{p_1 s_1} \theta_{s_1}^{q_1}, \\
\nabla c^{p_2}_{q_2i} = d c^{p_2}_{q_2i} - c^{p_2}_{q_2j}
\omega_i^j + c^{s_2}_{q_2i} \omega_{s_2}^{p_2}
- c^{p_2}_{s_2i} \theta^{s_2}_{q_2}, \\
\nabla b^\alpha_{p_2 q_2} = d b^\alpha_{p_2 q_2} + b^\beta_{s_2 q_2}
\omega_\beta^\alpha - b^\alpha_{s_2 q_2} \omega^{s_2}_{p_2}
- b^\alpha_{p_2 s_2} \theta^{s_2}_{q_2}.
\end{array}
\right.
\renewcommand{\arraystretch}{1}
$$

Consider the system of equations
\begin{equation}\label{eq:42}
   \theta^{q_1} = 0
\end{equation}
on the manifold $X$. By (8), its exterior differentiation gives
\begin{equation}\label{eq:43}
   \theta^{q_2} \wedge \theta_{q_2}^{q_1} = 0.
\end{equation}
It follows from (43) that the conditions of complete integrability
of equations (42) have the form
\begin{equation}\label{eq:44}
    \theta_{q_2}^{q_1} = l_{q_2 s_2}^{q_1} \theta^{s_2},
\end{equation}
where $ l_{q_2 s_2}^{q_1} = l_{s_2 q_2}^{q_1}.$

By equations (42), the system of equations (38) takes the form
\begin{equation}\label{eq:45}
 (c^{s_2}_{q_2i}   \omega^{p_1}_{s_2} - c^{p_1}_{s_1 i} \theta_{q_2}^{s_1})
 \wedge  \theta^{q_2} = 0.
\end{equation}
Suppose that the component $F_1$ of the focus hypersurface $F_L$
does not have multiple components. Assuming that $l \geq 1$, we
write equations (45) for two different values of the index $i$,
for example, for $i = 0, 1$. Since the matrices
$(c^{p_1}_{s_1i})$ and $(c^{p_2}_{s_2i})$ are not proportional,
then it follows from (45) that two terms occurring in (45)
vanish separately. In particular, this means that
\begin{equation}\label{eq:46}
c^{p_1}_{s_1 i} \theta_{q_2}^{s_1}
 \wedge  \theta^{q_2} = 0.
\end{equation}
Since the number of linearly independent forms
among the 1-forms $\omega_i^{p_1}$ connected with
the basis forms by relations (47) is equal to
the number of linearly independent forms $\theta^{q_1}$
(i.e., it is equal $r_1$), then it follows from (45) that
$$
\theta_{q_2}^{s_1}
 \wedge  \theta^{q_2} = 0.
$$
But the last equations coincide with equations (43) and
are conditions of complete integrability of (42). Thus
the manifold $X$ is foliated into an $r_1$-parameter family of
submanifolds of dimension $l + r_2$ and of rank $r_2$, and
these submanifolds belong to the types described in Theorems 2 or
3.

In a similar way, one can prove the complete
integrability of equations $\theta^{q_2} = 0$ on the manifold $X$.
Thus the manifold $X$ is foliated  also into an $r_2$-parameter family of
submanifolds of dimension $l + r_1$ and of rank $r_1$.

By induction over $s$, we can prove the result, which we have proved for
$s = 2$ components, for the case of any number $s$ of components.
\rule{3mm}{3mm}

Thus, Theorem 5 describes the structure of tangentially
degenerate manifolds of general types. As a result, this theorem
is a {\em structure theorem} for such manifolds.

Note that the torsal manifolds described in Theorem 2 are
completely reducible, and the manifolds $X$
described in Theorems 3 and 4 are irreducible manifolds.

Note that Theorem 5 does not cover  tangentially
degenerate submanifolds with multiple nonlinear components
of their focal images.
 This gives rise to  the following problem.

\textbf{Problem}. {\em Construct an example of a submanifold
$X \subset P^N ({\bf C})$ with
a degenerate Gauss map whose focal images have multiple nonlinear
components or prove that such submanifolds do not exist.}

\section{Additional results}

In conclusion we  prove two additional theorems.

\begin{theorem} Let $X \subset P^N$ be a tangentially
degenerate submanifold of dimension $n$ and rank $r < n$.
Suppose that
all matrices $B^\alpha$ can be simultaneously diagonalized,
$ B^\alpha = \mbox{{\rm diag}} \;
(b^\alpha_{l+1, l+1}, \ldots , b^\alpha_{nn})$.
Suppose also that the rectangular matrix $B$
$($defined by $(20))$ composed from the eigenvalues of
the matrices $B^\alpha$ has a rank $r_1 \leq r - 1$,
and this rank does not reduce when we delete
any column of this matrix. Then the submanifold $X$
belongs to a subspace $P^{n+r_1}$ of  the space $P^N$.
\end{theorem}

{\sf Proof.} Under the  conditions  of Theorem 11,
the second group of equations (11) takes the form
\begin{equation}\label{eq:47}
\omega_p^\alpha = b_{pp}^\alpha \theta^p, \;\;
p = l + 1, \ldots, n, \; \alpha = n + 1, \ldots , N.
\end{equation}
The matrix $B$ has only $r_1$ linearly
independent rows. Thus by means of transformations
of moving frame's vertices located outside
of the tangent subspace $T_L$, equations (47)
can be reduced to the form
\begin{equation}\label{eq:48}
\omega_p^\lambda = b_{pp}^\lambda \theta^p, \;\; \omega_p^\sigma = 0,
\end{equation}
where $\lambda = n + 1, \ldots, n + r_1, \;
\sigma = n + r_1 + 1, \ldots , N.$
The second group of equations (9) takes the form
$$
dA_p =  \omega_p^i A_i + \omega_p^q A_q +
 \omega_p^\lambda A_\lambda,
$$
and the points $A_\lambda$ together with the points
$A_i$ and $A_q$ define the osculating subspace $S_L$ of
the submanifold $X$ for all points $x \in L$.
The dimension of  $S_L$ is $n + r_1, \;
\dim \, S_L = n + r_1$.

Differentiation of the points $A_\lambda$ gives
\begin{equation}\label{eq:49}
dA_\lambda =  \omega_\lambda^i A_i + \omega_\lambda^p A_p +
 \omega_\lambda^\mu A_\mu + \omega_\lambda^\sigma A_\sigma,
\end{equation}
where $\lambda, \mu = n+1, \ldots , n + r_1; \;
\sigma = n + r_1 + 1, \ldots , N$. If $\theta^p = 0$,
then the osculating subspace $S_L$ of $X$ remains fixed.
It follows from equations (49) that the 1-forms $\omega_\lambda^\rho$ are
expressed in terms of the basis forms $\theta^p$
of $X$, that is,
\begin{equation}\label{eq:50}
\omega_\lambda^\rho = l_{\lambda p}^\rho  \theta^p.
\end{equation}

Taking exterior derivatives of the second group
of equations (48), we find that
\begin{equation}\label{eq:51}
\omega_p^\lambda \wedge \omega_\lambda^\rho = 0.
\end{equation}
Substituting the values of the 1-forms
$\omega_p^\lambda$ and $\omega_\lambda^\rho$ from equations
(48) and (50) into equation (51), we find that
$$
b_{pp}^\lambda \theta^p \wedge l_{\lambda q}^\rho \theta^q = 0.
$$
In this equation the summation is carried over
the indices $\lambda$ and $q$, but there is no summation over
the index $p$.
It follows from these equations that
\begin{equation}\label{eq:52}
b_{pp}^\lambda l_{\lambda q}^\rho = 0, \;\; p \neq q.
\end{equation}
System (52) is a system of linear homogeneous system with respect
to the unknown variables $ l_{\lambda q}^\rho$. For each pair
of the values $\rho$ and $q$,  system  (52) has the rank
$r - 1$ and $r_1$ unknowns. Since $r_1 \leq r - 1$,
under the conditions of Theorem 11, the rank
of the matrix of coefficients of this system is equal
$r_1$. As a result, the system has only the trivial solution
$ l_{\lambda q}^\rho = 0$. Thus equations (50) take  the form
\begin{equation}\label{eq:53}
\omega_p^\lambda = 0.
\end{equation}
It follows from (49) and (53) that the osculating subspace
$S_L$ of $X$ remains fixed when $L$ moves in $X$.
Thus $X \subset P^{n + r_1}$.
\rule{3mm}{3mm}

\textbf{Remark.}
If $r_1 = r$ and $N > n + r$, then the osculating subspace $S_L$ of $X$
can move in $P^N$ when $L$ moves in $X$. In this case the
submanifold $X$ is torsal.

Theorem 11 is similar to Theorem 3.10 from the book [AG 93] and
was proved in [AG 93] for submanifolds of a space $P^N$ bearing
a net of conjugate lines. Note that Theorem 3.10 from
[AG 93] generalizes a similar theorem of
C. Segre (see [Se 07], p. 571) proved for
 submanifolds $X$ of dimension $n$ of the space $P^N$
 which has at each point $x \in X$ the osculating
 subspace $S_x$ of dimension $n + 1$. By this theorem,
a submanifold $X$  either belongs to a subspace $P^{n+1}$
or is a torse.

The theorem dual to Theorem 11 is also valid.
\begin{theorem} Let $X \subset P^N$ be a tangentially
degenerate submanifold of dimension $n$ and rank $r < n$.
Suppose that
all  matrices $C_i$ can be simultaneously diagonalized,
 $ C_i = \mbox{{\rm diag}} \; (c^{l+1}_{l+1, i}, \ldots , c^n_{ni})$.
 Suppose also that the rectangular matrix $C = (c_{pi}^p)$
 composed from the eigenvalues of
the matrices $C_i$ has a rank $r_2 \leq r - 1$,
and this rank is not reduced when we delete
any column of this matrix. Then the submanifold $X$
is a cone with an $(l - r_2)$-dimensional vertex $K_L$.
\end{theorem}

{\sf Proof.} The proof of this theorem is similar to the
proof of Theorem 11. \rule{3mm}{3mm}

\noindent {\em Authors' addresses}:\\

\noindent
\begin{tabular}{ll}
M. A. Akivis &V. V. Goldberg\\
Department of Mathematics &Department of Mathematical Sciences\\
Jerusalem College of Technology---Mahon Lev &  New
Jersey Institute of Technology \\
Havaad Haleumi St., P. O. B. 16031 & University Heights \\
 Jerusalem 91160, Israel &  Newark, N.J. 07102, U.S.A. \\
 & \\
 E-mail address: akivis@avoda.jct.ac.il & E-mail address:
 vlgold@m.njit.edu
 \end{tabular}

\begin{thebibliography}{AG 93}

\bibitem[A 57]{A:[A 57]}
Akivis, M. A., {\em Focal images of a surface of rank $r$,} (Russian) Izv.
      Vyssh. Uchebn. Zaved. Mat. \textbf{1957}, no. 1, 9--19.

\bibitem[A 62]{A:[A 62]}
Akivis, M. A., {\em On a class of tangentially degenerate surfaces,}
      (Russian)
      Dokl. Akad. Nauk SSSR \textbf{146} (1962), no. 3, 515--518.
      English transl: Soviet Math. Dokl.
      {\bf 3} (1962), no. 5, 1328--1331.

\bibitem[A 87]{A:[A 87]}      Akivis, M. A.,
 {\em On multidimensional strongly parabolic surfaces,}
      (Russian) Izv. Vyssh.
       Uchebn. Zaved. Mat. \textbf{1987}, no. 5 (311), 3--10.
       English transl: Soviet Math. (Iz. VUZ) {\bf 31}
       (1987), no. 5, 1--11.

\bibitem[AG 93]{AG:[AG 93]}      Akivis, M. A., and V. V. Goldberg,
{\em Projective differential geometry and its generalizations},
North-Holland, Amsterdam, 1993, xii+362 pp.


\bibitem[AG 00]{AG:[AG 00]}      Akivis, M. A., and V. V. Goldberg,
{\em Equivalence of examples of Sacksteder and Bourgain},
7 pp., 2000 (submitted).


\bibitem[AGL]{AGL:[AGL]}      Akivis, M. A.,  V. V.
Goldberg, and J. Landsberg, {\em On the
Griffiths-Harris conjecture on varieties with degenerate Gauss
mappings}, 1999, 3 pp. (submitted).

\bibitem[AR 64]{AR:[AR 64]} Akivis, M. A. and  V. V. Ryzhkov,
 {\em Multidimensional surfaces of
     special projective types}, (Russian) Proc. Fourth
     All-Union Math. Congr. (Leningrad, 1961),
      Vol. II, pp. 159--164, Izdat. "Nauka",
     Leningrad, 1964.

 \bibitem[B 97]{Bo:[Bo 97]} Borisenko, A. A.,
 {\em Extrinsic geometry of strongly parabolic multidimensional
 submanifolds} (Russian), Uspekhi Mat.
Nauk \textbf{52} (1997), no. 6(318), 3--52; English transl:
Russian Math. Surveys \textbf{52} (1997), no. 6, 1141--1190.

\bibitem[Br 38]{B:[Br 38]}  Brauner, K.,
{\em \"{U}ber Mannigfaltigkeiten, deren
Tangentialmannigfaltigkeiten ausgeart sind}, Monatsh. Math.
Phys. \textbf{46} (1938), 335--365.


\bibitem[C 16]{C:[C 16]}
Cartan, \'{E}.,
{\em La deformation des hypersurfaces dans l'espace
euclidien reel \`{a} $n$ dimensions},   Bull. Soc.
Math. France \textbf{44} (1916), 65--99.

\bibitem[C 19]{C:[C 19]}
Cartan, \'{E}.,
 {\em Sur les vari\'{e}t\'{e}s de courbure
constante d'un espace euclidien ou non-euclidien},  Bull. Soc.
Math. France \textbf{47} (1919), 125--160; {\bf 48} (1920),
132--208.

\bibitem[C 45]{C:[C 45]}
Cartan, \'{E}., {\em Les syst\`{e}mes diff\'{e}rentiels ext\'{e}rieurs et
                  leurs applications g\'{e}om\'{e}triques},
                  Hermann, Paris, 1951, 214 pp.



\bibitem[CK 52]{CK:[CK 52]}
Chern, S. S. and  N. H. Kuiper,
  {\em Some theorems on isometric imbeddings of
compact Riemannian manifolds in Euclidean space},
Ann. of Math. (2) \textbf{56} (1952), 422--430.

\bibitem[D 89]{D:[D 89]}
Delano\"{e}, Ph.,
{\em L'op\'{e}rateur de Monge-Amp\`{e}re r\'{e}el
et la g\'{e}om\'{e}trie des sous-varo\'{e}t\'{e}s},
in {\em Geometry and Topology of Submanifolds} (Eds. J. M. Morvan
and L. Verstraelen), World Scientific, 1989, pp. 49--72.


\bibitem[DFN 85]{DFN:[DFN 85]}
Dubrovin, B. A., A. T. Fomenko, and S. P. Novikov, {\em Modern
geometry -- methods and applications, Part II. The
geometry and topology of
manifolds}, Springer-Verlag, New York--Berlin, 1985, xv+430 pp.

\bibitem[E 86]{E:[E 86]}Ein, L.,
{\em Varieties with small dual varieties. $I$},
Invent. Math. \textbf{86} (1986), no. 1, 63--74.

\bibitem[E 85]{E:[E 85]}Ein, L.,
{\em Varieties with small dual varieties. $II$},
Duke Math. J. \textbf{52} (1985), no. 4, 895--907.


\bibitem[FW 95]{FW:[FW 95]} Fischer, G. and H. Wu,
  {\em Developable complex analytic submanifolds},
 Internat. J. Math. \textbf{6} (1995), no. 2, 229--272.

\bibitem[GH 79]{GH:[GH 79]} Griffiths, P. A. and J. Harris,
{\em Algebraic geometry and local
differential geometry}, Ann. Sci. \'{E}cole Norm. Sup. (4)
\textbf{12} (1979), 355--452.

\bibitem[I 98]{I:[I 98]} Ishikawa, G.,
{\em Developable hypersurfaces and algebraic homogeneous spaces
in real projective space}, in
{\em Homogeneous structures and theory of submanifolds} (Kyoto, 1998).
S\=urikaisekikenky\=usho K\=oky\=uroku No. 1069 ,
(1998), 92--104.

\bibitem[I 99a]{I:[I 99a]} Ishikawa, G.,
{\em Singularities of developable surfaces}, Singularity
theory (Liverpool, 1996), xxii--xxiii, 403--418, London Math. Soc.
Lecture Note Ser., 263, Cambridge Univ. Press, Cambridge, 1999.

\bibitem[I 99b]{I:[I 99b]} Ishikawa, G.,
{\em Developable hypersurfaces and homogeneous spaces
in a real projective space}, Lobachevskii J. Math.
 \textbf{3} (1999), 113--125.


\bibitem[IM 97]{IM:[IM 97]} Ishikawa, G. and T. Marimoto,
{\em Solution surfaces of Monge-Amp\`{e}re equations}, Hokkaido
Univ. Preprint Series \textbf{376} (1997), 16 pp.

\bibitem[KN 69]{KN:[KN 69]}
Kobayashi, S. and  K. Nomizu, {\em Foundations of differential
geometry},  Vol. 2,  Wiley--Interscience, New York/London/Sydney,
1969, xv+470 pp.

\bibitem[L 96]{L:[L 96]}     Landsberg, J. M.,
{\em  On degenerate secant and tangential varieties and local
differential  geometry},   Duke Math. J. \textbf{85} (1996),
605--634.

\bibitem[L 99]{L:[L 99]}     Landsberg, J. M., {\em Algebraic
geometry and projective differential geometry},
Lecture Notes Series, No. 45, Seoul National Univ.,
Seoul, Korea, 1999, 85pp.

\bibitem[R 58]{R:[R 58]}
Ryzhkov, V. V.,
{\em Conjugate nets on multidimensional surfaces}, (Russian)
Trudy Moskov. Mat. Obshch. \textbf{7} (1958), 179--226.

\bibitem[R 60]{R:[R 60]}
 Ryzhkov, V. V.,
{\em Tangential degenerate surfaces}, (Russian) Dokl.
       Akad. Nauk SSSR \textbf{135} (1960), no. 1, 20--22.
       English  transl: Soviet Math. Dokl. \textbf{1} (1960), no. 1,
1233--1236.


\bibitem[S 60]{S:[S 60]}
Sacksteder, R: {\em On hypersurfaces with no negative sectional
curvature}, Amer. J. Math.  \textbf{82} (1960), no. 3, 609--630.


\bibitem[Sa 57]{Sa:[Sa 57]}
 Savelyev, S. I.,
{\em Surfaces with plane generators along which the tangent
plane is constant}, (Russian) Dokl.
       Akad. Nauk SSSR \textbf{115} (1957), no. 4, 663--665.

\bibitem[Sa 60]{Sa:[Sa 60]}
Savelyev, S. I.,
{\em On surfaces with plane generators along which the tangent
plane is constant}, (Russian)
 Izv. Vyssh. Uchebn. Zaved. Mat. \textbf{1960}, no. 2
(15), 154--167.

\bibitem[Se 07]{Se:[Se 07]}
Segre, C.,
{\em Su una classe di superficie degli iperspazi legate colle
equazioni lineari alle derivate parziali di $2^\circ$ ordine},
Atti Accad. Sci. Torino Cl. Sci. Fis. Mat. Natur.
 \textbf{42} (1907), 559--591.

\bibitem[W 77]{W:[W 77]} Wolf, J. A.,
{\em Spaces of constant curvature}, 4th ed., Publish or Perish,
Berkeley,  1977, xvi+408 pp.

\bibitem[Wu 95]{Wu:[Wu 95]} Wu, H.,
 {\em Complete developable submanifolds in real and complex Euclidean
 spaces}, Internat. J. Math. \textbf{6} (1995), no. 3, 461--489.


\bibitem[WZ 99]{WZ:[WZ 99]} Wu, H. and P. Zheng,
 {\em On complete developable submanifolds in complex Euclidean
 spaces}, Preprint, 30 pp. (July 15, 1999).

\bibitem[Y 53]{Y:[Y 53]} Yanenko, N. N.,
{\em Some questions of the theory of embeddings of Riemannian
metrics into Euclidean spaces}, (Russian) Uspekhi Mat. Nauk
\textbf{8} (1953), no. 1 (53), 21--100.

\end{thebibliography}
\end{document}